\documentclass[pdflatex,sn-mathphys-num]{sn-jnl}

\usepackage{graphicx}%
\usepackage{multirow}%
\usepackage{amsmath,amssymb,amsfonts}%
\usepackage{amsthm}%
\usepackage{mathrsfs}%
\usepackage{upgreek}
\usepackage[title]{appendix}%
\usepackage{xcolor}%
\usepackage{textcomp}%
\usepackage{manyfoot}%
\usepackage{booktabs}%
\usepackage{algorithm}%
\usepackage{algorithmicx}%
\usepackage{algpseudocode}%
\usepackage{listings}%

\theoremstyle{thmstyleone}%
\theoremstyle{thmstyletwo}%

\theoremstyle{thmstylethree}%

\begin{document}

\title[Article Title]{Ultra-Granular Calorimeter Performances for the Heavy Flavor Physics Program at the Z Peak }

\author{J.-C.~Brient}\email{brient@llr.in2p3.fr}

\affil*{\orgdiv{Laboratoire Leprince-Ringuet}, \orgname{CNRS-IPP}, \orgaddress{\street{Route de Saclay}, \city{Palaiseau}, \postcode{91128},  \country{France}}}


\abstract{ Various decays of the B mesons are here used to establish the performances of an ultra-granular electromagnetic calorimeter for heavy flavour physics at an electron positron accelerator running at the Z peak. The silicon-tungsten electromagnetic calorimeter of the ILD concept is used for this purpose, enhanced by a timing capability. 
 When possible, a $\pi^0$ mass fit of $\gamma \gamma$ system is performed to improve the $\pi^0$ energy resolution. It is also shown that in the presence in the final state of a photon without a $\pi^0$, the identification of genuine photon(s) versus fake photon(s) coming from  K$_\mathrm{L}$'s, neutron's, debris of hadronic shower or other high energy $\pi^0$, is essential. It allows for a good precision and a good signal over background ratio for this kind of physics. The possible impact for the tau physics is discussed. 
}


%

\keywords{Ultra-granular ECAL, heavy flavour physics}

\maketitle

\section{Introduction}\label{sec1}

The next generation of electron-positron collider will be able to produce a wealth of events when running at the Z peak. Such a high luminosity would give access to a large programme of heavy flavour physics with an excellent precision~\cite{bartmann_future_2025}. The detector requirements for this physics have for long been based on the paradigm that a very good electromagnetic calorimeter energy resolution is needed\cite{bib1}.  Several concepts of detectors and technologies have been proposed for the future colliders, among which imaging electromagnetic calorimeters (ECAL), with rather modest energy resolution. However, imaging calorimeters can provide excellent precisions on the position of photon showers. Based on a silicon-tungsten calorimeter~\cite{boudry_new_2023} proposed for ILD~\cite{abramowicz_ild_2025} and the recent improvement of the time measurement for each pixel~\cite{acar_timing_2024, HGTD}, a new photon reconstruction algorithm~\cite {bib2} has been developed, based on causality between pixels. This algorithm is very efficient at very low-energy, while also providing a very high purity, i.e.\ no fake “photon” coming from all primary sources ($\mathrm{K_L}$, neutron, debris of the hadronic shower from charged particles). It also has the capability of analyzing close-by showers.

For the study presented here, a fast simulation has been developed based on the performances of the ILD concept and of this algorithm. When $\pi^0$(s) are involved, this study will try to find the optimal photon pairing and do a 1C-fit of the $\pi^0$ mass for the $\gamma \gamma$ system. 
 
By selecting a few processes from the list of benchmark processes suggested for the heavy flavour physics on these leptonic colliders, the study will measure the improvement associated with the new algorithm's good precision on the position provided when added to the 1C-fit of the $\pi^0$ mass. 

For all these analyses, the signal and backgrounds have been generated with PYTHIA~\cite{sjostrand_introduction_2015}.In the following, a reference electromagnetic calorimeter energy resolution of $\Delta E/E=5\,\%/\sqrt{E}$, is used for the comparisons as a representative benchmark.

\section{Improvement due to the $\pi^0$ mass fit}
\label{pi0_mass_fit}

\subsection{ Preliminary test with isolated $\pi^0$'s} 
A preliminary test has been done  using single isolated $\pi^0$s at low energy (below $10\,$GeV). Taking the 2 photons and performing a mass fit of the $\gamma \gamma$ system, using the expected performances on the photon position for the silicon-tungsten ECAL, lead to an impressive improvement.
 The energy resolution of  $\Delta E/E=3\,\%/\sqrt{E}$, quoted in the figure  corresponds only to the performance assumed in the CEPC study ~\cite{bib3}, and is used here for qualitative comparison.

Our first conclusion is clear: when possible, the use of  $\pi^0$  mass fit, leads to a large improvement on the energy resolution of the  $\pi^0$, with a result comparable to the resolution related in~\cite{bib3}.

\subsection{Bs and B$^0$ decaying to Ds$\pi$ or DsK} 
\begin{figure}
    \centering
    \includegraphics[width=0.9\linewidth]{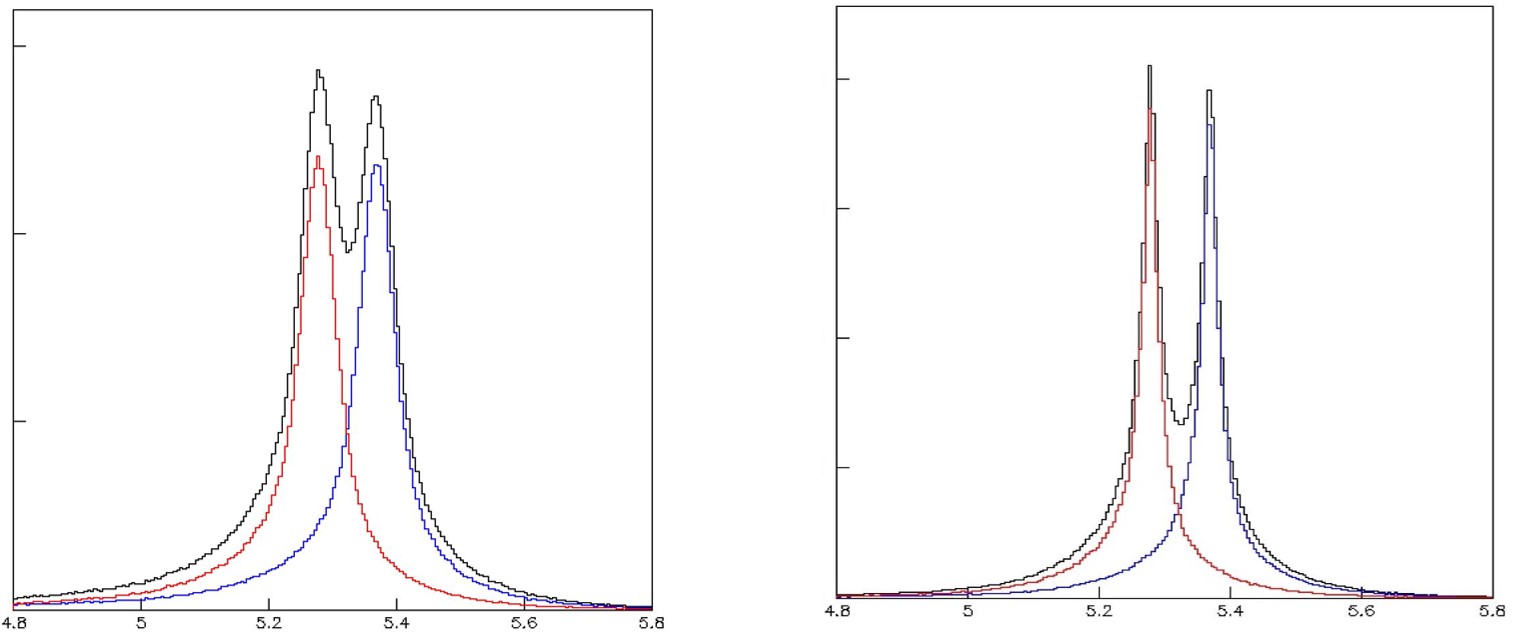}
    \caption{The Ds $\pi$ mass distribution with Bs and B$^0$ decays. On the left when using an ECAL with an energy resolution of $\Delta E/E= 5\%/\sqrt{E}$. On the right, with an ultra-granular calorimeter with $\Delta E/E= 16\%/\sqrt{E} $  but with an additional 1C-fit of the $\pi^0$ mass }
    \label{fig:bsbo}
\end{figure}

For this study, the Ds decay to $\phi\,\rho$ is used, with $\phi$ decaying to K$^+$K$^-$  and $\rho$ to $\pi^\pm \pi^0$. For this decays chain, the Ds $\pi$  gives a final state with 2 charged kaons and 2 charged pions, together with 2 photons coming from the $\pi^0$. For the analysis, the events are split into 2 hemispheres using the thrust axis. Working on the hemispheres separately, a simple and crude selection  is performed to isolate the 2 “good” photons from a large number  coming from the signal but also from fragmentation. After this selection, the comparison of the $\phi \rho$ mass spectra for different calorimeters shows that the  resolution of the $\phi \rho$ mass is particularly bad for the raw silicon-tungsten ECAL, but also for a crystal calorimeter with 3\% resolution. Only the fit of the $\pi^0$  mass gives a  distribution good enough to select the Ds signal. Of course, the impact on the separation Bs /B$^0$ using the mass of the system is also significant. It is illustrated on the Figure~\ref{fig:bsbo}. 
Defining the standard deviation of the mass resolution spectrum as the quadratic sum of the mass resolution for Bs and for B$^0$, we can estimate the separation Bs/B$^0$ at $0.6\,\sigma$ for the calorimeter with $\Delta E/E=5\,\%/\sqrt{E}$ and as $1.6\,\sigma$ when a $\pi^0$ mass fit is performed. 
The best separation between Bs and B$^0$ is clearly obtained with an ultra-granular ECAL with a position precision such that a $\pi^0$ mass fit can be performed. It is not discussed here but, for this study, the shower energy threshold to identify a photon could be a major problem for many proposed ECAL. This is not a problem for ultra-granular ECAL, since the new reconstruction algorithm  is very efficient, with a perfect purity (no fake photon is observed using GEANT4 simulation and the new photon reconstruction algorithm used in this study), down to 50 MeV
\cite {bib2}.

\subsection{  B$^0$ decays to $\pi^0\pi^0$     } 
\begin{figure}
        \centering
        \includegraphics[width=0.5\linewidth]{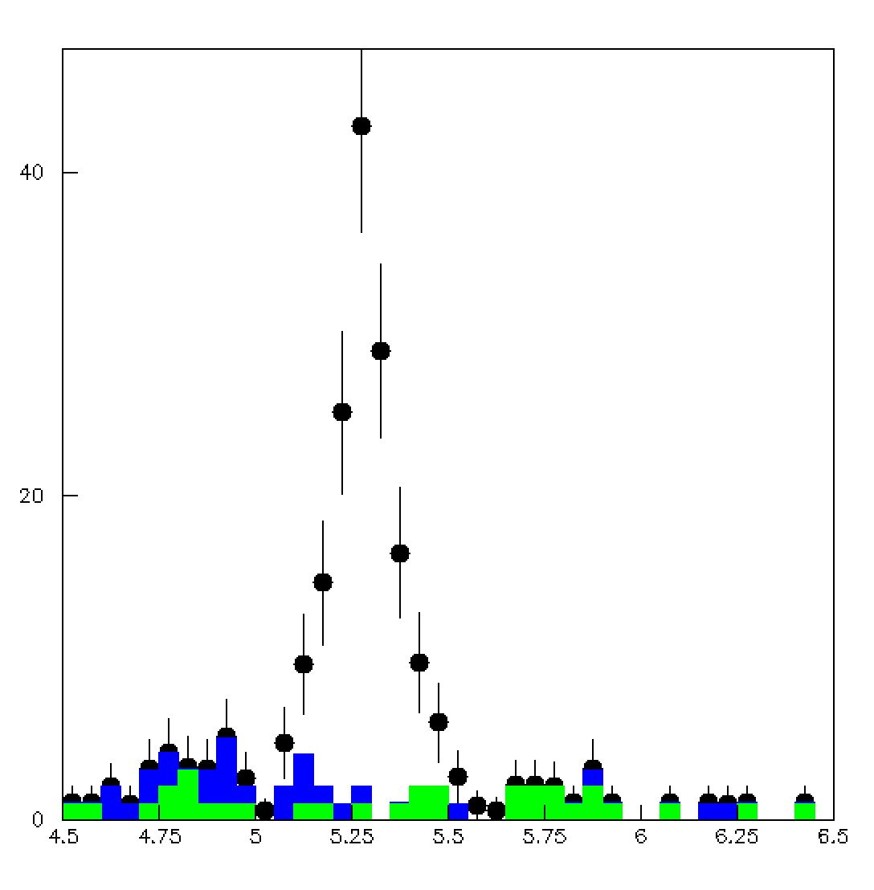}
        \caption{Based on $10^9$ Z decays, the figure shows the distribution of the 2 $\pi^0$  mass, after the 1C-fit of each $\pi^0$ (see text). The blue histogram contains the overall background from b and c quark Z decays, the green one the charm background. The black dots are the expected signal.}
        \label{fig:b2p}
\end{figure}

For this process, the events are split in 2 hemispheres, using the thrust axis. Working on each hemisphere separately, one is tagged as the “b hemisphere” and in the opposite, “the search hemisphere”, we require no vertex with charged tracks to be reconstructed at large distance from the interaction point. The relevant parameters for the photons to be signed are the position precision and the separability between a single photon and a high-energy $\pi^0$.
Taking the highest energy photon in the “search hemisphere” and its closest  “companion”, a good signal of $\pi^0$  is observed on the $\gamma\gamma$  mass spectrum. Among the remaining photons of the hemisphere, the one with the largest energy is selected and then, pairing it with the  other remaining photons of the same hemisphere, we select those pairs with a mass close to the $\pi^0$ mass, together with an overall 4 $\gamma$'s system with a mass between 4 and 6 GeV. This allows for a good selection of the signal. For each $\gamma\gamma$ system, a 1C-fit to the $\pi^0$ mass is then performed.  Having tagged the b quark on the opposite hemisphere, we can readily neglect the u, d, s background. An efficiency of $90\,\%$ on the b-tagging and a $10\,\%$ c quark misidentification are assumed. Simple cuts are applied on top of those described (for example, on the energy of the $\pi^0\pi^0$ system).  Using BR(B$^0 \to \pi^0\pi^0) = 1.55\times 10^{-6}$ as quoted in the PDG, we obtain the mass distribution shown in Figure~\ref{fig:b2p}. An excellent signal over background is possible, as soon as each 
 $\pi^0$ mass is fitted on each of the selected $\gamma$$\gamma$ system. 
 
\section{Improvements related to the $\gamma$$/\pi^0$ identification}\label{sec2}

\subsection{b$\to$ s$\gamma$ transition study,  using  B$^0$ $\to \gamma$K* } 
For this study, a 1C-fit of any $\pi^0$ mass is not possible, but the size of the pixels in the ultra-granular ECAL ($0.5\times0.5\,\mathrm{cm^2}$ silicon pixel for ILD), allows for an excellent $\gamma/\pi^0$ separation, easily up to $\pi^0$'s energies of $35\,$GeV.    With a more sophisticated identification such as Bulos shower moments ~\cite{bulos} , based on longitudinal and transverse shower shape variables, allowing discrimination between single photons and overlapping photons from high-energy  $\pi^0$  decays, it could be possible disentangle single photon from  $\pi^0$ , at even at larger $\pi^0$  energy

\begin{figure}
        \centering
        \includegraphics[width=0.5\linewidth]{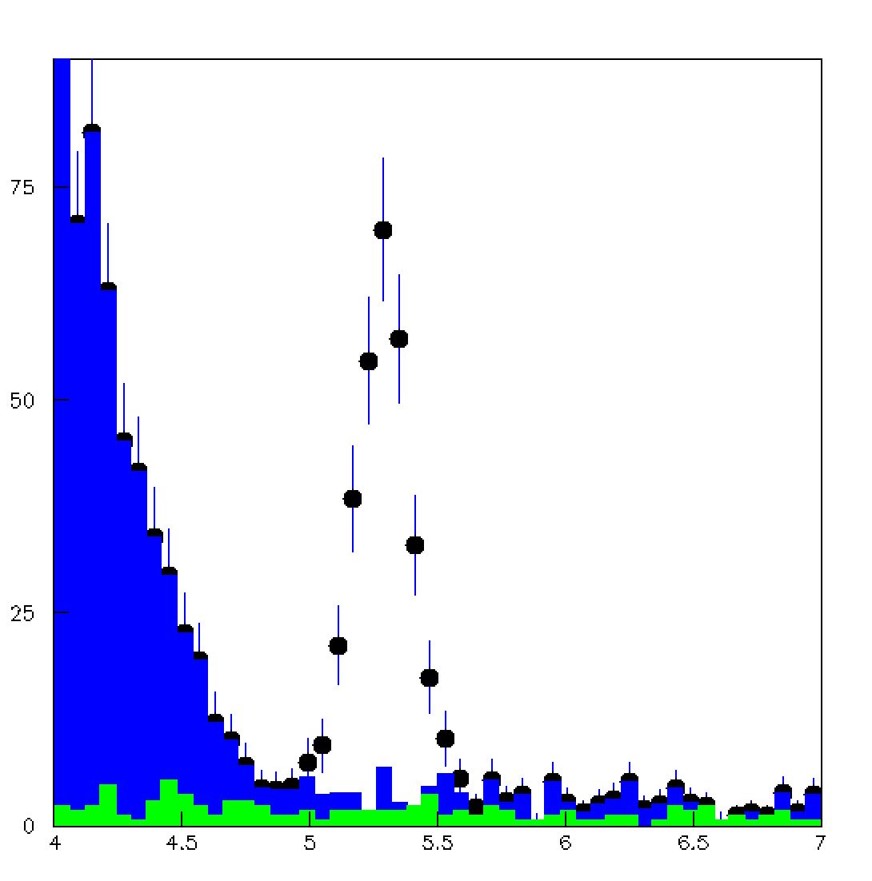}
        \caption{Distributions of the (K$^\pm$ pi$^\mp$ $\gamma$) mass for 10$^9$ Z decays mass. In blue the total b and c background, the green histogram is the c background and the black dots are the expected signal. }
    \label{fig:kpgamma}
\end{figure}

This separation is essential since there are many background decays of B meson with a displaced vertex from K$^*$(890) and at least one $\pi^0$.
For the selection, a displaced vertex (K$^\pm$ $\pi^\mp$) is selected and the correct high-energy photon is identified in the jet cone.
Then, some cuts are applied, first, for the  displaced vertex (K$^\pm\pi^\mp$), requiring a distance larger than 50$\,\upmu$m from the interaction point.
Then, the  (K$^\pm$ $\pi^\mp$ ) mass has to be in the K$^*$ range,  within $0.85–1.0\,$GeV. Finally, the photon energy is required to be larger than $5\,$GeV and  the energy of the (K$^*\gamma$) system larger than $30\,$GeV. The expected mass distributions for $10^9$ Z decays are shown in Figure~\ref{fig:kpgamma}.

\section{Tau physics }\label{sec12}
Concerning rare decays of the $\tau$ lepton, the possibility to fit the mass of a $\pi^0$ or a $\eta$ could help a lot, but the purity of the photon sample, together with a high reconstruction efficiency, are the most important parameters. 
From very accurate experiments done in the past to measure tau polarisation~\cite{heister_measurement_2001}, we know that the most important factors are the efficiency and purity of the photon(s) reconstructed samples at both very low and very high $\pi^0$ energies. This is done to get rid of sources of electromagnetic showers that aren't caused by photons.
The new photon reconstruction algorithm seems perfect for that purpose. However, further studies have to be done to confirm the cross talk between tau decay channels and the precision on the variable relevant for measuring the tau polarisation in the channels $\tau \to \rho\nu$ and $\tau \to \pi\nu$.

\section{Conclusion}\label{sec13}
From preliminary analyses on various heavy flavour final states, for a collider running at the Z-peak centre-of-mass energy, it can be concluded that an ultra-granular electromagnetic calorimeter -- such as the one proposed to ILD, based on silicon pixels and tungsten radiator, with a timing resolution comparable to the one of ATLAS HGTD -- offers superior performances. 
This is especially visible in the final state of the B meson decays with at least one $\pi^0$, when a 1C-Fit of the mass of the $\pi^0$ is performed. In channels where the relevant parameter is the purity of the photon sample and especially when the signal/background is strongly related to the separation of single $\gamma$ from high energy $\pi^0$'s, the potential of high-granularity with timing capacity is overwhelming. The possible impact for tau physics is also discussed. 

\backmatter





\bmhead{Acknowledgements}

We acknowledge Henri Videau, in particular, for the new photon reconstruction algorithm, and Vincent Boudry for constructive
discussions.


\bibliography{sn-bibliography}


\begin{thebibliography}{11}
\ifx \bisbn   \undefined \def \bisbn  #1{ISBN #1}\fi
\ifx \binits  \undefined \def \binits#1{#1}\fi
\ifx \bauthor  \undefined \def \bauthor#1{#1}\fi
\ifx \batitle  \undefined \def \batitle#1{#1}\fi
\ifx \bjtitle  \undefined \def \bjtitle#1{#1}\fi
\ifx \bvolume  \undefined \def \bvolume#1{\textbf{#1}}\fi
\ifx \byear  \undefined \def \byear#1{#1}\fi
\ifx \bissue  \undefined \def \bissue#1{#1}\fi
\ifx \bfpage  \undefined \def \bfpage#1{#1}\fi
\ifx \blpage  \undefined \def \blpage #1{#1}\fi
\ifx \burl  \undefined \def \burl#1{\textsf{#1}}\fi
\ifx \doiurl  \undefined \def \doiurl#1{\url{https://doi.org/#1}}\fi
\ifx \betal  \undefined \def \betal{\textit{et al.}}\fi
\ifx \binstitute  \undefined \def \binstitute#1{#1}\fi
\ifx \binstitutionaled  \undefined \def \binstitutionaled#1{#1}\fi
\ifx \bctitle  \undefined \def \bctitle#1{#1}\fi
\ifx \beditor  \undefined \def \beditor#1{#1}\fi
\ifx \bpublisher  \undefined \def \bpublisher#1{#1}\fi
\ifx \bbtitle  \undefined \def \bbtitle#1{#1}\fi
\ifx \bedition  \undefined \def \bedition#1{#1}\fi
\ifx \bseriesno  \undefined \def \bseriesno#1{#1}\fi
\ifx \blocation  \undefined \def \blocation#1{#1}\fi
\ifx \bsertitle  \undefined \def \bsertitle#1{#1}\fi
\ifx \bsnm \undefined \def \bsnm#1{#1}\fi
\ifx \bsuffix \undefined \def \bsuffix#1{#1}\fi
\ifx \bparticle \undefined \def \bparticle#1{#1}\fi
\ifx \barticle \undefined \def \barticle#1{#1}\fi
\bibcommenthead
\ifx \bconfdate \undefined \def \bconfdate #1{#1}\fi
\ifx \botherref \undefined \def \botherref #1{#1}\fi
\ifx \url \undefined \def \url#1{\textsf{#1}}\fi
\ifx \bchapter \undefined \def \bchapter#1{#1}\fi
\ifx \bbook \undefined \def \bbook#1{#1}\fi
\ifx \bcomment \undefined \def \bcomment#1{#1}\fi
\ifx \oauthor \undefined \def \oauthor#1{#1}\fi
\ifx \citeauthoryear \undefined \def \citeauthoryear#1{#1}\fi
\ifx \endbibitem  \undefined \def \endbibitem {}\fi
\ifx \bconflocation  \undefined \def \bconflocation#1{#1}\fi
\ifx \arxivurl  \undefined \def \arxivurl#1{\textsf{#1}}\fi
\csname PreBibitemsHook\endcsname

\bibitem[\protect\citeauthoryear{Bartmann et~al.}{2025}]{bartmann_future_2025}
\begin{botherref}
\oauthor{\bsnm{Bartmann}, \binits{W.}}, et al.:
Future {{Circular Collider Feasibility Study Report Volume}} 1: {{Physics}} and
  {{Experiments}}.
Technical report,
CERN Document Server
(2025).
\doiurl{10.17181/CERN.9DKX.TDH9} .
\url{http://cds.cern.ch/record/2928193}
\end{botherref}
\endbibitem

\bibitem[\protect\citeauthoryear{Dams}{Nov. 2025}]{bib1}
\begin{botherref}
\oauthor{\bsnm{Dams}, \binits{M.}}:
Detector challenges and concepts at \uppercase{FCC}-ee
(Nov. 2025).
presented at Flavours Physics at \uppercase{FCC} Workshop at CERN
\end{botherref}
\endbibitem

\bibitem[\protect\citeauthoryear{Boudry}{2023}]{boudry_new_2023}
\begin{barticle}
\bauthor{\bsnm{Boudry}, \binits{V.}}:
\batitle{New results of the technological prototype of the {{CALICE}} highly
  granular silicon tungsten calorimeter}.
\bjtitle{Nucl. Instrum. Meth. A}
\bvolume{1051},
\bfpage{168185}
(\byear{2023})
\doiurl{10.1016/j.nima.2023.168185}
\end{barticle}
\endbibitem

\bibitem[\protect\citeauthoryear{Abramowicz et~al.}{2025}]{abramowicz_ild_2025}
\begin{botherref}
\oauthor{\bsnm{Abramowicz}, \binits{H.}}, et al.:
The {{ILD Detector}}: {{A Versatile Detector}} for an {{Electron-Positron
  Collider}} at {{Energies}} up to 1 {{TeV}}.
arXiv
(2025).
\doiurl{10.48550/arXiv.2506.06030} .
\url{http://arxiv.org/abs/2506.06030}
Accessed 2025-11-19
\end{botherref}
\endbibitem

\bibitem[\protect\citeauthoryear{Acar et~al.}{2024}]{acar_timing_2024}
\begin{barticle}
\bauthor{\bsnm{Acar}, \binits{B.}},
\bauthor{\bsnm{others}},
\bauthor{\bsnm{{collaboration}}, \binits{T.C.H.}}:
\batitle{Timing performance of the {{CMS High Granularity Calorimeter}}
  prototype}.
\bjtitle{Journal of Instrumentation}
\bvolume{19}(\bissue{04}),
\bfpage{04015}
(\byear{2024})
\doiurl{10.1088/1748-0221/19/04/P04015} .
Accessed 2026-01-20
\end{barticle}
\endbibitem

\bibitem[\protect\citeauthoryear{Missio}{2024}]{HGTD}
\begin{barticle}
\bauthor{\bsnm{Missio}, \binits{M.}}:
\batitle{{Overview of the ATLAS High-Granularity Timing Detector: project
  status and results}}.
\bjtitle{JINST}
\bvolume{19}(\bissue{04}),
\bfpage{04008}
(\byear{2024})
\doiurl{10.1088/1748-0221/19/04/C04008}
\end{barticle}
\endbibitem

\bibitem[\protect\citeauthoryear{Videau}{2025}]{bib2}
\begin{botherref}
\oauthor{\bsnm{Videau}, \binits{H.}}:
New photon(s) reconstruction algorithm based on ultragranular {ECAL} with a
  good timing precision
(2025).
Private communication, to be published.
\end{botherref}
\endbibitem

\bibitem[\protect\citeauthoryear{Sj{\"o}strand
  et~al.}{2015}]{sjostrand_introduction_2015}
\begin{barticle}
\bauthor{\bsnm{Sj{\"o}strand}, \binits{T.}}, \betal:
\batitle{An introduction to {{PYTHIA}} 8.2}.
\bjtitle{Comput. Phys. Commun.}
\bvolume{191},
\bfpage{159}--\blpage{177}
(\byear{2015})
\doiurl{10.1016/j.cpc.2015.01.024}
\end{barticle}
\endbibitem

\bibitem[\protect\citeauthoryear{Guo}{2025}]{bib3}
\begin{botherref}
\oauthor{\bsnm{Guo}, \binits{F.}}:
Sim-digi-rec software activities for \uppercase{CEPC} calorimeters
(2025).
DRD 6 meeting at Orsay, April. 2025
\end{botherref}
\endbibitem

\bibitem[\protect\citeauthoryear{F.Bossi et~al.}{1993}]{bulos}
\begin{botherref}
\oauthor{\bsnm{F.Bossi}}, et al.:
Minute of the photon reconstruction meeting.
Aleph 91-133,
CERN,
CERN
(1993).
\url{https://cds.cern.ch/record/805407/files/aleph-91-133.pdf}
\end{botherref}
\endbibitem

\bibitem[\protect\citeauthoryear{Heister
  et~al.}{2001}]{heister_measurement_2001}
\begin{barticle}
\bauthor{\bsnm{Heister}, \binits{A.}}, \betal:
\batitle{Measurement of the tau polarization at {{LEP}}}.
\bjtitle{Eur. Phys. J. C}
\bvolume{20},
\bfpage{401}--\blpage{430}
(\byear{2001})
\doiurl{10.1007/s100520100689}
\end{barticle}
\endbibitem

\end{thebibliography}

\end{document}